\documentclass[12pt,reqno]{amsart}
\usepackage{amssymb}
\title[Families with a strictly maximal Higgs field]{Families over curves
with a strictly maximal Higgs field}
\author[Eckart Viehweg]{Eckart Viehweg}
\address{Universit\"at Essen, FB6 Mathematik, 45117 Essen, Germany}
\email{ viehweg@uni-essen.de}
\thanks{This work has been supported by the ``DFG-Schwerpunktprogramm
Globale Methoden in der Komplexen Geometrie''. The second named author
is supported by a grant from the Research
Grants Council of the Hong Kong
Special Administrative Region, China
(Project No. CUHK 4034/02P)}
\author[Kang Zuo]{Kang Zuo}
\address{The Chinese University of Hong Kong, Department of Mathematics,
Shatin, Hong Kong}
\email{kzuo@math.cuhk.edu.hk}
\begin{document}
\theoremstyle{plain}
\newtheorem{thm}{Theorem}[section]
\newtheorem{theorem}[thm]{Theorem}
\newtheorem{lemma}[thm]{Lemma}
\newtheorem{corollary}[thm]{Corollary}
\newtheorem{proposition}[thm]{Proposition}
\theoremstyle{definition}
\newtheorem{notations}[thm]{Notations}
\newtheorem{problem}[thm]{Problem}
\newtheorem{remark}[thm]{Remark}
\newtheorem{remarks}[thm]{Remarks}
\newtheorem{definition}[thm]{Definition}
\newtheorem{claim}[thm]{Claim}
\newtheorem{assumption}[thm]{Assumption}
\newtheorem{assumptions}[thm]{Assumptions}
\newtheorem{properties}[thm]{Properties}
\newtheorem{example}[thm]{Example}
\numberwithin{equation}{thm}
\catcode`\@=11
\def\opn#1#2{\def#1{\mathop{\kern0pt\fam0#2}\nolimits}}
\def\bold#1{{\bf #1}}%
\def\underrightarrow{\mathpalette\underrightarrow@}
\def\underrightarrow@#1#2{\vtop{\ialign{$##$\cr
 \hfil#1#2\hfil\cr\noalign{\nointerlineskip}%
 #1{-}\mkern-6mu\cleaders\hbox{$#1\mkern-2mu{-}\mkern-2mu$}\hfill
 \mkern-6mu{\to}\cr}}}
\let\underarrow\underrightarrow
\def\underleftarrow{\mathpalette\underleftarrow@}
\def\underleftarrow@#1#2{\vtop{\ialign{$##$\cr
 \hfil#1#2\hfil\cr\noalign{\nointerlineskip}#1{\leftarrow}\mkern-6mu
 \cleaders\hbox{$#1\mkern-2mu{-}\mkern-2mu$}\hfill
 \mkern-6mu{-}\cr}}}
\let\amp@rs@nd@\relax
\newdimen\ex@
\ex@.2326ex
\newdimen\bigaw@
\newdimen\minaw@
\minaw@16.08739\ex@
\newdimen\minCDaw@
\minCDaw@2.5pc
\newif\ifCD@
\def\minCDarrowwidth#1{\minCDaw@#1}
\newenvironment{CD}{\@CD}{\@endCD}
\def\@CD{\def\A##1A##2A{\llap{$\vcenter{\hbox
 {$\scriptstyle##1$}}$}\Big\uparrow\rlap{$\vcenter{\hbox{%
$\scriptstyle##2$}}$}&&}%
\def\V##1V##2V{\llap{$\vcenter{\hbox
 {$\scriptstyle##1$}}$}\Big\downarrow\rlap{$\vcenter{\hbox{%
$\scriptstyle##2$}}$}&&}%
\def\={&\hskip.5em\mathrel
 {\vbox{\hrule width\minCDaw@\vskip3\ex@\hrule width
 \minCDaw@}}\hskip.5em&}%
\def\verteq{\Big\Vert&&}%
\def\noarr{&&}%
\def\vspace##1{\noalign{\vskip##1\relax}}\relax\let\amp@rs@nd@&\iffalse}\fi
 \CD@true\vcenter\bgroup\relax\let\\=\cr\iffalse}\fi\tabskip\z@skip\baselineskip20\ex@
 \lineskip3\ex@\lineskiplimit3\ex@\halign\bgroup
 &\hfill$\m@th##$\hfill\cr}
\def\@endCD{\cr\egroup\egroup}
\def\>#1>#2>{\amp@rs@nd@\setbox\z@\hbox{$\scriptstyle
 \;{#1}\;\;$}\setbox\@ne\hbox{$\scriptstyle\;{#2}\;\;$}\setbox\tw@
 \hbox{$#2$}\ifCD@
 \global\bigaw@\minCDaw@\else\global\bigaw@\minaw@\fi
 \ifdim\wd\z@>\bigaw@\global\bigaw@\wd\z@\fi
 \ifdim\wd\@ne>\bigaw@\global\bigaw@\wd\@ne\fi
 \ifCD@\hskip.5em\fi
 \ifdim\wd\tw@>\z@
 \mathrel{\mathop{\hbox to\bigaw@{\rightarrowfill}}\limits^{#1}_{#2}}\else
 \mathrel{\mathop{\hbox to\bigaw@{\rightarrowfill}}\limits^{#1}}\fi
 \ifCD@\hskip.5em\fi\amp@rs@nd@}
\def\<#1<#2<{\amp@rs@nd@\setbox\z@\hbox{$\scriptstyle
 \;\;{#1}\;$}\setbox\@ne\hbox{$\scriptstyle\;\;{#2}\;$}\setbox\tw@
 \hbox{$#2$}\ifCD@
 \global\bigaw@\minCDaw@\else\global\bigaw@\minaw@\fi
 \ifdim\wd\z@>\bigaw@\global\bigaw@\wd\z@\fi
 \ifdim\wd\@ne>\bigaw@\global\bigaw@\wd\@ne\fi
 \ifCD@\hskip.5em\fi
 \ifdim\wd\tw@>\z@
 \mathrel{\mathop{\hbox to\bigaw@{\leftarrowfill}}\limits^{#1}_{#2}}\else
 \mathrel{\mathop{\hbox to\bigaw@{\leftarrowfill}}\limits^{#1}}\fi
 \ifCD@\hskip.5em\fi\amp@rs@nd@}
\newenvironment{CDS}{\@CDS}{\@endCDS}
\def\@CDS{\def\A##1A##2A{\llap{$\vcenter{\hbox
 {$\scriptstyle##1$}}$}\Big\uparrow\rlap{$\vcenter{\hbox{%
$\scriptstyle##2$}}$}&}%
\def\V##1V##2V{\llap{$\vcenter{\hbox
 {$\scriptstyle##1$}}$}\Big\downarrow\rlap{$\vcenter{\hbox{%
$\scriptstyle##2$}}$}&}%
\def\={&\hskip.5em\mathrel
 {\vbox{\hrule width\minCDaw@\vskip3\ex@\hrule width
 \minCDaw@}}\hskip.5em&}
\def\verteq{\Big\Vert&}
\def\novarr{&}
\def\noharr{&&}
\def\SE##1E##2E{\slantedarrow(0,18)(4,-3){##1}{##2}&}
\def\SW##1W##2W{\slantedarrow(24,18)(-4,-3){##1}{##2}&}
\def\NE##1E##2E{\slantedarrow(0,0)(4,3){##1}{##2}&}
\def\NW##1W##2W{\slantedarrow(24,0)(-4,3){##1}{##2}&}
\def\slantedarrow(##1)(##2)##3##4{%
\thinlines\unitlength1pt\lower 6.5pt\hbox{\begin{picture}(24,18)%
\put(##1){\vector(##2){24}}%
\put(0,8){$\scriptstyle##3$}%
\put(20,8){$\scriptstyle##4$}%
\end{picture}}}
\def\vspace##1{\noalign{\vskip##1\relax}}\relax\let\amp@rs@nd@&\iffalse}\fi
 \CD@true\vcenter\bgroup\relax\let\\=\cr\iffalse}\fi\tabskip\z@skip\baselineskip20\ex@
 \lineskip3\ex@\lineskiplimit3\ex@\halign\bgroup
 &\hfill$\m@th##$\hfill\cr}
\def\@endCDS{\cr\egroup\egroup}
\newdimen\TriCDarrw@
\newif\ifTriV@
\newenvironment{TriCDV}{\@TriCDV}{\@endTriCD}
\newenvironment{TriCDA}{\@TriCDA}{\@endTriCD}
\def\@TriCDV{\TriV@true\def\TriCDpos@{6}\@TriCD}
\def\@TriCDA{\TriV@false\def\TriCDpos@{10}\@TriCD}
\def\@TriCD#1#2#3#4#5#6{%
\setbox0\hbox{$\ifTriV@#6\else#1\fi$}
\TriCDarrw@=\wd0 \advance\TriCDarrw@ 24pt
\advance\TriCDarrw@ -1em
\def\SE##1E##2E{\slantedarrow(0,18)(2,-3){##1}{##2}&}
\def\SW##1W##2W{\slantedarrow(12,18)(-2,-3){##1}{##2}&}
\def\NE##1E##2E{\slantedarrow(0,0)(2,3){##1}{##2}&}
\def\NW##1W##2W{\slantedarrow(12,0)(-2,3){##1}{##2}&}
\def\slantedarrow(##1)(##2)##3##4{\thinlines\unitlength1pt
\lower 6.5pt\hbox{\begin{picture}(12,18)%
\put(##1){\vector(##2){12}}%
\put(-4,\TriCDpos@){$\scriptstyle##3$}%
\put(12,\TriCDpos@){$\scriptstyle##4$}%
\end{picture}}}
\def\={\mathrel {\vbox{\hrule
   width\TriCDarrw@\vskip3\ex@\hrule width
   \TriCDarrw@}}}
\def\>##1>>{\setbox\z@\hbox{$\scriptstyle
 \;{##1}\;\;$}\global\bigaw@\TriCDarrw@
 \ifdim\wd\z@>\bigaw@\global\bigaw@\wd\z@\fi
 \hskip.5em
 \mathrel{\mathop{\hbox to \TriCDarrw@
{\rightarrowfill}}\limits^{##1}}
 \hskip.5em}
\def\<##1<<{\setbox\z@\hbox{$\scriptstyle
 \;{##1}\;\;$}\global\bigaw@\TriCDarrw@
 \ifdim\wd\z@>\bigaw@\global\bigaw@\wd\z@\fi
 \mathrel{\mathop{\hbox to\bigaw@{\leftarrowfill}}\limits^{##1}}
 }
 \CD@true\vcenter\bgroup\relax\let\\=\cr\iffalse}\fi
 \tabskip\z@skip\baselineskip20\ex@
 \lineskip3\ex@\lineskiplimit3\ex@
 \ifTriV@
 \halign\bgroup
 &\hfill$\m@th##$\hfill\cr
#1&\multispan3\hfill$#2$\hfill&#3\\
&#4&#5\\
&&#6\cr\egroup%
\else
 \halign\bgroup
 &\hfill$\m@th##$\hfill\cr
&&#1\\%
&#2&#3\\
#4&\multispan3\hfill$#5$\hfill&#6\cr\egroup
\fi}
\def\@endTriCD{\egroup}
\newcommand{\sA}{{\mathcal A}}
\newcommand{\sB}{{\mathcal B}}
\newcommand{\sC}{{\mathcal C}}
\newcommand{\sD}{{\mathcal D}}
\newcommand{\sE}{{\mathcal E}}
\newcommand{\sF}{{\mathcal F}}
\newcommand{\sG}{{\mathcal G}}
\newcommand{\sH}{{\mathcal H}}
\newcommand{\sI}{{\mathcal I}}
\newcommand{\sJ}{{\mathcal J}}
\newcommand{\sK}{{\mathcal K}}
\newcommand{\sL}{{\mathcal L}}
\newcommand{\sM}{{\mathcal M}}
\newcommand{\sN}{{\mathcal N}}
\newcommand{\sO}{{\mathcal O}}
\newcommand{\sP}{{\mathcal P}}
\newcommand{\sQ}{{\mathcal Q}}
\newcommand{\sR}{{\mathcal R}}
\newcommand{\sS}{{\mathcal S}}
\newcommand{\sT}{{\mathcal T}}
\newcommand{\sU}{{\mathcal U}}
\newcommand{\sV}{{\mathcal V}}
\newcommand{\sW}{{\mathcal W}}
\newcommand{\sX}{{\mathcal X}}
\newcommand{\sY}{{\mathcal Y}}
\newcommand{\sZ}{{\mathcal Z}}
\newcommand{\A}{{\mathbb A}}
\newcommand{\B}{{\mathbb B}}
\newcommand{\C}{{\mathbb C}}
\newcommand{\D}{{\mathbb D}}
\newcommand{\E}{{\mathbb E}}
\newcommand{\F}{{\mathbb F}}
\newcommand{\G}{{\mathbb G}}
\newcommand{\HH}{{\mathbb H}}
\newcommand{\I}{{\mathbb I}}
\newcommand{\J}{{\mathbb J}}
\newcommand{\LL}{{\mathbb L}}
\newcommand{\M}{{\mathbb M}}
\newcommand{\N}{{\mathbb N}}
\newcommand{\BP}{{\mathbb P}}
\newcommand{\Q}{{\mathbb Q}}
\newcommand{\R}{{\mathbb R}}
\newcommand{\Shur}{{\mathbb S}}
\newcommand{\T}{{\mathbb T}}
\newcommand{\U}{{\mathbb U}}
\newcommand{\V}{{\mathbb V}}
\newcommand{\W}{{\mathbb W}}
\newcommand{\X}{{\mathbb X}}
\newcommand{\Y}{{\mathbb Y}}
\newcommand{\Z}{{\mathbb Z}}
\newcommand{\id}{{\rm id}}
\newcommand{\rank}{{\rm rank}}
\newcommand{\END}{{\mathbb E}{\rm nd}}
\newcommand{\End}{{\rm End}}
\newcommand{\Hg}{{\rm Hg}}
\newcommand{\rk}{{\rm rank}}
\newcommand{\Sl}{{\rm Sl}}
\newcommand{\Gl}{{\rm Gl}}
\newcommand{\Cor}{{\rm Cor}}
\newcommand{\width}{{\rm width}}
\begin{abstract}
We study certain rigid Shimura curves in the moduli scheme of
polarized minimal $n$-folds of Kodaira dimension zero. Those are characterized
by some numerical condition on the Deligne extension of the corresponding
variation of Hodge structures, or equivalently by the strict maximality
of the induced Higgs field. We show that such Shimura curves can not be proper
for $n$ odd, and we give some examples, showing that they exist in all dimensions.
\end{abstract}
\maketitle
Let $f:X\to Y$ be a family of $n$-dimensional complex algebraic varieties,
smooth over $U=Y\setminus S$. We will assume that $X$ is projective
and non-singular, that $Y$ is a smooth projective curve, and that the
general fibre $F$ of $f$ is connected, hence irreducible. Writing
$X_0=f^{-1}(U)$ one has the $\Q$ variations of Hodge
structures $R^kf_*\Q_{X_0}$ of weight $k$. Usually we will assume that the
monodromy around each point $s\in S$ is unipotent and that the family
is not birationally isotrivial.

Considered the Deligne extension of $(R^kf_*\Q_{X_0})\otimes \sO_U$ to $Y$ together
with the extension of the Hodge filtration. Taking the graded sheaf
one obtains the Higgs bundle
$$
(E,\theta)=(\bigoplus_{p+q=k} E^{p,q},\theta^{p,q})
$$
with $E^{p,q}=R^q f_*\Omega^p_{X/Y}(\log \Delta)$, for $\Delta=f^*(S)$.
The restriction $\theta^{q,p}$ of the Higgs field to $E^{q,p}$
is given by the edge morphisms
$$
R^q f_*\Omega^p_{X/Y}(\log \Delta)\>>> R^{q+1}f_*\Omega^{p-1}_{X/Y}  \otimes
\Omega^{1}_{Y}(\log S)
$$
of the $p$-th wedge product of the tautological sequence
$$
0\to {f}^*\Omega^1_Y(\log S)\to \Omega^1_{X}(\log \Delta) \to
\Omega_{X/Y}^1(\log \Delta))\to 0.
$$
In \cite{VZ4} we studied families of Abelian varieties, and $k=1$.
Then $E^{1,0}$ is a direct sum $F^{1,0} \oplus N^{1,0}$
with $F^{1,0}$ ample and $N^{1,0}$ flat, hence $N^{1,0}\subset{\rm
Ker}(\theta^{1,0})$. Correspondingly, one has
$E^{0,1}=F^{0,1}\oplus N^{0,1}$ and $E$ is the direct sum of the Higgs bundles
$$
(F=F^{1,0}\oplus F^{0,1}, \theta^{1,0}|_{F^{1,0}}) \mbox{ \ \ \
and \ \ \ } (N^{1,0}\oplus N^{0,1},0).
$$
In this special case we called the Higgs field maximal if
$$
\theta^{1,0}:F^{1,0} \to F^{0,1}\otimes \Omega^1_Y(\log S)
$$
is an isomorphism. The notion ``strictly maximal'', introduced below in \ref{maxhiggs},
adds the additional condition that $N^{1,0}=N^{0,1}=0$, or in different terms that
$R^1f_*\Q_{X_0}$ has no unitary part.

The Arakelov inequalities (\cite{Fal}, generalized in \cite{De4}, \cite{Pet2},
\cite{JZ}) say that
$$
2\cdot \deg(E^{1,0}) \leq {\rm rank}(E^{1,0}) \cdot \deg(\Omega^1_Y(\log S)),
$$
and this inequality becomes an equality if and only if the
Higgs field is strictly maximal.

In this article we extend both, the Arakelov inequalities, and the
concept of the strict maximality to variation of Hodge structures $\X$
of weight $k\geq 1$, and we discuss some implications for the shape
of the variation of Hodge structures, for their special Mumford-Tate groups
and for the rigidity of the family. Let us remark already, that
for $k$ even one allows $R^kf_*\Q_{X_0}$ to have a unitary part $\U$, as long as
it is concentrated in bidegree $(\frac{k}{2},\frac{k}{2})$. Our main interest
are variations of Hodge structures of weight $k=n$, given by a family
$f:X\to Y$ with general fibre $F$ a minimal model. We doubt that there are
examples with $\kappa(F)>0$. For $\kappa(F)=0$ a non-isotrivial family
with a strictly maximal Higgs field gives rise to a Shimura curve in the moduli space
of polarized manifolds, and parallel to the case of Abelian varieties one is tempted
to expect that there are few of such curves.

After giving the definition of strictly maximal Higgs fields
in Section \ref{secmax} we will prove in Section \ref{secarak} a generalized
Arakelov inequality for variations of Hodge structures of weight $k\geq 0$.
As for $k=1$, these inequalities become equalities if and only if the
Higgs field is strictly maximal.

C. Simpson's theory of Higgs bundles will allow in Section \ref{secdecomp}
to construct certain decompositions for Higgs bundles with a strictly maximal
Higgs field. Roughly speaking, as for Abelian varieties \cite{VZ4} or for
$K3$-surfaces \cite{STZ}, they are decomposed in direct sums and tensor products of
a fixed $\C$-variation of Hodge structures $\LL$ of weight one
and some unitary parts. This decomposition will imply rigidity and some minimality
for the special Mumford-Tate group in Section \ref{sectrigi}.

In case $S\neq \emptyset$ the unipotence of the local monodromies
will imply in Section \ref{secsplit} that both, the local system
$\LL$ and the decomposition can be defined over $\Q$. This
implies that the variation of Hodge structures looks like
one given by the $n$-th cohomology of the $n$-th product
$E\times_Y \cdots \times_Y E$, where $E\to Y$ is a modular family
of elliptic curves. For $S=\emptyset$ and for $k$ odd, we will show
next that there are no variations of Hodge structures with a strictly
maximal Higgs field and with $h^{k,0}=1$.

In the final Section \ref{secex} we will give examples of families
$f:X\to Y$ of $n$-dimensional varieties, whose $n$-th variation of Hodge structures
is strictly maximal. Those families are obtained as quotients of families
of Abelian varieties, and their general fibre $F$ is a minimal model of
Kodaira dimension zero. For $S=\emptyset$ we only have examples of families
of Kummer type. In particular they only exist for $n$ even, and they
are not Calabi-Yau manifolds, except for $n=2$.
For $S\neq \emptyset$, generalizing a construction due to
C. Borcea \cite{Bor}, will construct families of Calabi-Yau varieties
by taking certain quotients of products of modular families of elliptic curves,
as predicted by Theorem \ref{nonempty}.\\

This note grew out of discussions started when the first named author visited
the Morningside Center of Mathematics in the Chinese Academy of Sciences in
Beijing. The final version was written during his visit of the Institute of
Mathematical Science and the Department of Mathematics at the Chinese University
of Hong Kong. He would like to thank the members of those Institutes for their
hospitality and help.

\section{Strictly maximal Higgs fields}
\label{secmax}
\begin{definition}\label{gmaxhiggs}
Let $\X$ be a polarized $\C$ variation of Hodge structures of
weight $k$, with Higgs bundle
$$(E,\theta)=(\bigoplus_{p+q=k} E^{p,q},\theta^{p,q}).
$$
$\X$ (or $(E,\theta)$) has a strict generically maximal Higgs field,
if one has a direct sum decomposition
\begin{equation}\label{dec1}
(E,\theta)=\bigoplus_{i=0}^k
(F_i,\theta|_{F_i}),
\end{equation}
with:
\begin{enumerate}
\item[i.] If $k-i$ is odd, $F_i = 0$.
\item[ii.] If $k-i$ is even, $F_i^{p,k-p}=F_i\cap E^{p,k-p}=0$
for $2p < k-i$ and $2p > k + i$.
Moreover, for all $p$ with $k-i+2 \leq 2p \leq k + i$,
$$
\tau_i^{p,k-p}=\theta|_{F_i^{p,k-p}}:F_i^{p,k-p}
\>>> F_i^{p-1,k-p+1}\otimes
\Omega_Y^1(\log S)
$$
is generically an isomorphisms.
\end{enumerate}
\end{definition}
So for $k-i$ even, $F_i^{p,k-p}$ can only be non-zero for
$$
p= \ \ \frac{k+i}{2} \ \ \Rightarrow \frac{k+i-2}{2} \ \ \Rightarrow \ \
\cdots \ \ \Rightarrow \ \ \frac{k-i+2}{2} \ \ \Rightarrow \ \ \frac{k-i}{2}
$$
and all Higgs fields \ \ \ $* \to *\otimes \Omega_Y^1(\log S)$ \ \ \
(indicated by $\Rightarrow$) are generically isomorphisms.
In particular, if $F_i\neq 0$ the width of $(F_i,\tau_i)$ is $i$.

The polarization of $\X$ induces an isomorphisms
$$
(F_i,\tau_i)
\simeq (F_i^\vee,\tau_i^\vee)(-k)
$$
where $(-k)$ stands for the Tate twist, i.e. for the shift of the
bigrading by $(k,k)$. In particular $\det(F_i)^2=\sO_Y$.

\begin{definition}\label{maxhiggs}
$\X$ (or $(E,\theta)$) has a strictly maximal Higgs field, if the morphisms
$\tau_i^{p,k-p}$ in \ref{gmaxhiggs}, ii), are all isomorphisms.
We will also say that $\X$ (or $(E,\theta)$) is strictly maximal.
\end{definition}
The strict maximality of the Higgs field allows to apply Simpson's Correspondence
\cite{Si2} (see also \cite{VZ4} and the references given there).
\begin{lemma}\label{deg}
Assume that $(E,\theta)$ has a strictly maximal Higgs field and let
$F_i$ be the maximal component of $E$ of width $i$. Then $\deg(F_i)=0$, and
$(F_i,\tau_i=\theta|_{F_i})$ is the Higgs bundle corresponding to a polarized
variation of Hodge structures $\V_i$.
\end{lemma}
\begin{proof} Write $\theta^{(p)}$ for the iterated map
$$
E \to E\otimes \Omega^1_Y(\log S) \to E\otimes \Omega^1_Y(\log S)^2 \to \cdots \to
E\otimes \Omega^1_Y(\log S)^p.
$$
For $i=k$ the sub Higgs bundle $(F_k,\tau_k)$ of $(E,\theta)$ is given by
$$F_k^{k-p,p}=\theta^{(p)}(E^{k,0})\otimes \Omega^1_Y(\log S)^{-p}.$$
Hence $\deg(F_k^{k-p,p})=\deg(E^{k,0}) - p\cdot \deg(\Omega^1_Y(\log
S))$ and for $p=k$ one finds
$$
-\deg(E^{k,0})=\deg(E^{0,k})=\deg(E^{k,0}) - k \cdot \deg(\Omega^1_Y(\log S).
$$
Adding up, one finds
\begin{gather*}
\deg(F_k)= k\cdot \deg(E^{k,0}) - \sum_{p=0}^{k}p\cdot \deg(\Omega^1_Y(\log
S))=\\
k\cdot \deg(E^{k,0}) - \frac{k(k+1)}{2} \cdot \deg(\Omega^1_Y(\log
S))= 0.
\end{gather*}
By Simpson's Correspondence one has a decomposition $\X=\V_k\oplus
\W$, which one may choose orthogonal with respect to the
polarization. $\W$ is again a variation of Hodge structures,
concentrated in bidegrees
$$(k-1,1) , \ldots , (1,k-1).$$
The Tate twist $\W(1)$ is a variation of Hodge structures of weight
$k-2$, obviously again with a strictly maximal Higgs field. By
induction on the weight of the variation of Hodge structures one
obtains \ref{deg}.
\end{proof}
\begin{definition}\label{maxpure} \
\begin{enumerate}
\item[a.] Assume that $\X$ has a strictly maximal Higgs field.
We call $(F_i,\tau_i)$ in \ref{gmaxhiggs} (or the local system $\V_i$
corresponding to the Higgs bundle $(F_i,\tau_i)$)
the strictly maximal pure component of $\X$ of width $i$.
\item[b.] $\X$ has a strictly maximal and pure Higgs field
(or $\X$ is a strictly maximal and pure variation of Hodge structures)
if it only has one strictly maximal pure component, i.e. if $\X=\V_i$, for some $i$.
\end{enumerate}
\end{definition}
So $\X$ is a strictly maximal and pure variation of Hodge structures,
if and only if all the
$$
\theta^{p,q}:E^{p,q} \>>> E^{p-1,q+1}\otimes \Omega_Y^1(\log S)
$$
are isomorphisms, except if $E^{p,q}=0$ or $E^{p-1,q+1}=0$.
\begin{remark}\label{maxhiggs2}
The decomposition (\ref{dec1}) gives rise to a
decomposition of $\C$ local systems:
\begin{gather}\label{dec2}
\X = \bigoplus_{i=0}^{[\frac{k-1}{2}]} \V_{2i+1}, \text{ \ \ for $k$ odd,}\\
\X = \bigoplus_{i=0}^{[\frac{k}{2}]} \V_{2i}, \text{ \ \ for $k$ even,}
\end{gather}
where the $\V_j$ are all zero or strictly maximal and pure of width $j$.
In particular, $\X$ can only have a unitary sub-system $\U$ for
$k$ even, and $\U=\V_0$ in this case.

By \cite{VZ4}, 3.3, in case $\X$ is a polarized $\bar\Q\cap \R$
variation of Hodge structures, the decomposition (\ref{dec2})
can be defined again over $\bar\Q\cap \R$, and orthogonal
with respect to the polarization.
\end{remark}
\begin{remarks}\label{stz}
As mentioned in the introduction, in \cite{VZ4} for the maximality of Higgs field,
of variation of Hodge structures of weight $1$,
one also allows unitary sub-systems of width $1$, and
in \cite{STZ} for families of K3-surfaces one also allows
Higgs fields with zero second iterated Kodaira Spencer map.
Those are again maximal in a certain sense, but not strictly maximal.
\end{remarks}

\section{Arakelov inequalities}
\label{secarak}
In this section we consider the Higgs bundle
$$
(E,\theta)=(\bigoplus_{p+q=k} E^{p,q},\theta^{p,q})
$$
corresponding to the Deligne extension of a polarized complex variation of
Hodge structures $\X$ of weight $k$ on $Y\setminus S.$
We will assume that the local monodromies around the points in $S$ are
all unipotent. Let us write
$$ E_0^{p,q}={\rm ker}(\theta^{p,q}: E^{p,q}\to E^{p-1,q+1}\otimes
\Omega^1_Y(\log S)),\quad h^{p,q}_0={\rm rk} E_0^{p,q}.$$
\cite{JZ} and \cite{Pet2} contains the proof of generalized Arakelov
inequalities:\\
\ \\
If $k=2l+1$, then
$$ \deg  E^{k,0}\leq \big(\frac{1}{2}(h^{k-l,l}-h_0^{k-l,l})+
\sum_{j=0}^{l-1}(h^{k-j,j}-h^{k-j,j}_0)\big)\cdot \deg(\Omega^1_Y(\log S)).$$
If $k=2l$,
$$ \deg E^{k,0}\leq \sum_{j=0}^{l-1}(h^{j,k-j}-h^{j,k-j}_0)\cdot
\deg(\Omega^1_Y(\log S)).$$

We will show a slightly different inequality, for which the upper bound is reached,
if and only if the Higgs field is strictly maximal.

\begin{proposition}\label{arakineq} \
\begin{enumerate}
\item[a.] For all $\nu \leq [\frac{k}{2}]$ one has
$$
\deg(E^{k-\nu,\nu}) \leq \frac{k-2\nu}{2}\cdot(h^{k-\nu,\nu}-h_0^{k-\nu,\nu})
\cdot \deg(\Omega^1_Y(\log S)).
$$
\item[b.] One has
$$
0\leq \sum_{\nu=0}^{[\frac{k}{2}]} \deg(E^{k-\nu,\nu}) \leq
\sum_{\nu=0}^{[\frac{k}{2}]}
\frac{k-2\nu}{2}\cdot(h^{k-\nu,\nu}-h_0^{k-\nu,\nu}) \cdot \deg(\Omega^1_Y(\log S)).
$$
\item[c.] Let $\mu$ be the smallest natural number with
$\theta^{k-\mu,\mu} \neq 0$. Then
$$
0 < \deg(E^{k-\mu,\mu}).
$$
In particular $\deg(\Omega^1_Y (\log S))>0$, if $\theta\neq 0$.\\
\item[d.] If $\theta\neq 0$ the following conditions are equivalent:\\
\begin{enumerate}
\item[(i)] For all $\nu \leq [\frac{k}{2}]$
$$
\deg(E^{k-\nu,\nu})= \frac{k-2\nu}{2}\cdot h^{k-\nu,\nu} \cdot \deg(\Omega^1_Y(\log S)).
$$
\item[(ii)]
$$
\sum_{\nu=0}^{[\frac{k}{2}]} \deg(E^{k-\nu,\nu})=
\sum_{\nu=0}^{[\frac{k}{2}]}
\frac{k-2\nu}{2}\cdot h^{k-\nu,\nu} \cdot \deg(\Omega^1_Y(\log S)).
$$
\item[(iii)] $(E,\theta)$ has a strictly maximal Higgs field.
\end{enumerate}
\end{enumerate}
\end{proposition}
If $f:X\to Y$ is a non-isotrivial family with smooth part
$V\to U$, and if $\omega_{F}$ is semi-ample for a general fibre $F$ of
$f$, then by \cite{VZ1}
$$\deg(\Omega^1_Y (\log S))>0.$$
\ref{arakineq}, c), is the corresponding statement for variations of Hodge
structures.
\begin{proof}[Proof of \ref{arakineq}]
Fix some $\nu \leq [\frac{k}{2}]$ and write
$$
A=\bigoplus_{q=0}^{\nu} E^{k-q,q}=E^{k,0}\oplus \cdots \oplus
E^{k-\nu,\nu}.
$$
The sub Higgs bundle $<A>$ of $E$, generated by $A$, contains $A$ as a
direct factor. Writing $<A>= A\oplus B'$, and
$\theta^{(j)}$ for the iterated map
$$
E \to E\otimes \Omega^1_Y(\log S) \to E\otimes \Omega^1_Y(\log S)^2 \to \cdots \to
E\otimes \Omega^1_Y(\log S)^j,
$$
one has
$$
B'= \bigoplus_{j=1}^{k-\nu}\theta^{(j)}(E^{k-\nu,\nu})\otimes \Omega^1_Y(\log S)^{-j}
= \bigoplus_{q=\nu+1}^{k}B^{k-q,q} \subset \bigoplus_{q=\nu+1}^{k} E^{k-q,q}.
$$
By definition $A\oplus B'$ is a sub Higgs bundle of $E$, hence \cite{Si2} implies
$$
\deg(A)+\deg(B') \leq 0.
$$
On the other hand,
$$
A^\vee=\bigoplus_{q=0}^{\nu} {E^{k-q,q}}^\vee = \bigoplus_{q=0}^{\nu} E^{q,k-q}
= \bigoplus_{q=k-\nu}^{k}E^{k-q,q}
$$
also is a Higgs sub-bundle, and
$$
-\deg(A^\vee)=\deg(A)=\sum_{q=0}^{\nu} {E^{k-q,q}} \geq 0.
$$
For $\nu=[\frac{k}{2}]$ we obtain the first inequality in
b). The second one is the sum of the inequalities
in a), which we will verify below.\\

For all $\nu \leq [\frac{k}{2}]$ \cite{Si2} implies that
$\deg(A^\vee)=0$, if and only if
$(A^\vee,\theta|_{A^\vee})$ is a direct factor of $(E,\theta)$.
Equivalently $(A,\theta|_{A})$ is a direct factor of $(E,\theta)$, or $B'=0$.
The latter is equivalent to $\theta^{(1)}(E^{k-\nu,\nu})=0$
hence to $\theta^{k-\nu,\nu}=0$.

In c) we assumed that $\theta^{k-\mu-i,\mu+i}=0$ for all $i>0$,
whereas $\theta^{k-\mu,\mu}\neq 0$. This implies
\begin{gather*}
\sum_{q=1}^{\mu-i}\deg(E^{k-q,q}) = 0, \mbox{ \ \ for \ } i>0
\text{ \ \ and \ \ } \sum_{q=1}^{\mu}\deg(E^{k-q,q}) > 0,
\end{gather*}
and thereby $\deg(E^{k-\mu,\mu})>0$.\\

Let us consider next the Higgs bundle
\begin{gather*}
A\oplus B \oplus \bigoplus_{q=k-\nu+1}^{k}E^{k-q,q}\text{ \ \ where}\\
B= \bigoplus_{q=\nu+1}^{k-\nu}B^{k-q,q}=
\bigoplus_{j=1}^{k-2\nu}\theta^{(j)}(E^{k-\nu,\nu})\otimes \Omega^1_Y(\log S)^{-j}.
\end{gather*}
Up to one factor $E^{\nu,k-\nu}$ the last term of this decomposition
coincides with $A^\vee$, and one finds
\begin{gather}
\deg(E^{k-\nu,\nu}) = \deg(A)+\sum_{q=k-\nu+1}^{k}E^{k-q,q} \leq -\deg(B) =\notag\\
\label{ara4}
- \sum_{q=\nu+1}^{k-\nu} \deg(B^{k-q,q})=
-\sum_{j=1}^{k-2\nu}\deg(B^{k-\nu-j,\nu+j}).
\end{gather}
By the choice of the $B^{q,p}$ the Higgs field $\theta$ induces surjections
\begin{gather}\label{surj}
E^{k-\nu,\nu} \>>> B^{k-\nu-1,\nu+1} \otimes \Omega^1_Y(\log S), \mbox{ \ \ and, for }
j > 0\\
B^{k-\nu-j,\nu+j} \>>> B^{k-\nu-j-1,\nu+j+1} \otimes \Omega^1_Y(\log S)\label{surj2}
\end{gather}
whose kernels, as a Higgs sub-bundles of $E$, have a non-positive degree.
Writing $B^{k-\nu,\nu}=E^{k-\nu,\nu}$, for a moment, one obtains for $j\geq 0$
\begin{multline*}
\deg(B^{k-\nu-j,\nu+j}) \leq\\
\deg(B^{k-\nu-j-1,\nu+j+1}) + \rk(B^{k-\nu-j-1,\nu+j+1})
\cdot \deg(\Omega^1_Y(\log S)),
\end{multline*}
and thereby
\begin{gather*}
\deg(E^{k-\nu,\nu}) \leq
\deg(B^{k-\nu-1,\nu+1}) + \rk(B^{k-\nu-1,\nu+1})\cdot \deg(\Omega^1_Y(\log S))\\
\leq \deg(B^{k-\nu-j,\nu+j}) +
\sum_{i=1}^{j} \rk(B^{k-\nu-i,\nu+i}) \cdot \deg(\Omega^1_Y(\log S)).
\end{gather*}
So the right hand side of (\ref{ara4}) has an upper bound
\begin{gather*}
-\sum_{j=1}^{k-2\nu}\deg(B^{k-\nu-j,\nu+j})
\leq\\
-(k-2\nu)\deg(E^{k-\nu,\nu}) +\sum_{j=1}^{k-2\nu}\sum_{i=1}^{j}
\rk(B^{k-\nu-i,\nu+i}) \cdot \deg(\Omega^1_Y(\log S)).
\end{gather*}
Composing the surjections (\ref{surj}) and (\ref{surj2}) one obtains
surjections
$$
E^{k-\nu,\nu} \>>> B^{k-\nu-i,\nu+i} \otimes \Omega^1_Y(\log S)^j.
$$
In particular $\rk(B^{k-\nu-i,\nu+i})$ is smaller than or equal to
$h^{k-\nu,\nu}-h_0^{k-\nu,\nu}$. Altogether one obtains
\begin{gather*}
\deg(E^{k-\nu,\nu}) \leq \hspace*{10cm}\\  -(k-2\nu)\deg(E^{k-\nu,\nu}) +
\sum_{j=1}^{k-2\nu} \sum_{i=1}^{j}\rk(B^{k-\nu-i,\nu+i}) \cdot
\deg(\Omega^1_Y(\log S))\leq \\
-(k-2\nu)\deg(E^{k-\nu,\nu}) + \sum_{j=1}^{k-2\nu}
\sum_{i=1}^{j}(h^{k-\nu,\nu}-h_0^{k-\nu,\nu}) \cdot \deg(\Omega^1_Y(\log S))=\\
-(k-2\nu)\deg(E^{k-\nu,\nu}) + \sum_{j=1}^{k-2\nu}
j\cdot(h^{k-\nu,\nu}-h_0^{k-\nu,\nu}) \cdot \deg(\Omega^1_Y(\log S))=\\
-(k-2\nu)\deg(E^{k-\nu,\nu}) +\hspace*{9cm}\\
\frac{(k-2\nu)(k-2\nu+1)}{2}(h^{k-\nu,\nu}-h_0^{k-\nu,\nu}) \cdot
\deg(\Omega^1_Y(\log S)),
\end{gather*}
and
$$
\deg(E^{k-\nu,\nu}) \leq \frac{(k-2\nu)}
{2}(h^{k-\nu,\nu}-h_0^{k-\nu,\nu}) \cdot \deg(\Omega^1_Y(\log S))
$$
as claimed in a).\\

If $\theta\neq 0$ one finds some $\mu$ satisfying the condition pose in c).
The inequality a) for $\nu=\mu$ implies that $\deg(\Omega^1_Y(\log S)) >0$.\\
Using a) and b) one finds that d), (i), and d), (ii), are equivalent,
and obviously both hold true if $(E,\theta)$ has a strictly maximal Higgs field.

Since we assume in d) that $\theta\neq 0$,
c) implies that $\deg(\Omega^1_Y(\log S)) >0$.
It remains to show that this condition together with the equations
in d), (i), force the Higgs field to be strictly maximal.

By a), for $\nu=0,\ldots,[\frac{k}{2}]$ one finds
\begin{gather}
\deg(E^{k-\nu,\nu}) \leq
\frac{k-2\nu}{2}\cdot(h^{k-\nu,\nu}-h_0^{k-\nu,\nu}) \cdot\deg(\Omega^1_Y(\log
S)) \leq \notag\\
\label{ara7} \frac{k-2\nu}{2}\cdot h^{k-\nu,\nu}
\cdot\deg(\Omega^1_Y(\log S)).
\end{gather}
The equations d), (i), say that both inequalities in
(\ref{ara7}) are equalities, hence $h_0^{k-\nu,\nu}=0$. So/
$$
\rk(B^{k-\nu-j-1,\nu+j+1})=h^{k-\nu,\nu},
$$
and the surjections (\ref{surj}) and (\ref{surj2})
are isomorphisms. Hence for all $\nu$
\begin{equation}\label{ara8}
\deg(E^{k-\nu,\nu}) = \deg(B^{k-\nu-j,\nu+j}) + j\cdot
h^{k-\nu,\nu}\cdot\deg(\Omega^1_Y(\log S)).
\end{equation}
As in the proof of \ref{deg} the equality (\ref{ara8}) allows recursively to
decompose the Higgs bundle $(E,\theta)$:\\

For $\nu=0$ choose the sub Higgs bundle
$(F_k,\tau_k)$ of with $F_k^{k,0}=E^{k,0}$ and with
$$
F_k^{k-j,j}=B^{k-j,j}=\theta^j(E^{k,0})\otimes
\Omega^1_Y(\log S)^{-j}.
$$
By (\ref{ara8}) and d), (i), for $\nu=0$,
\begin{gather*}
\deg(F_k)=\sum_{j=0}^{k} \deg(F_k^{k-j,j})= (k+1) \cdot \deg(E^{k,0}) -
\sum_{j=1}^{k} j\cdot \deg(\Omega^1_Y(\log S))=\\
(k+1) \cdot \deg(E^{k,0}) -
\frac{k(k+1)}{2} j\cdot \deg(\Omega^1_Y(\log S))=0.
\end{gather*}
By Simpson's Correspondence \cite{Si2} there exists a decomposition
$\X=\V_k\oplus \W$, where $\W$ is a variation of Hodge structures,
concentrated in bidegrees
$$
(k-1,1), \ldots , (1,k-1),
$$
and where $\V_k$ is strictly maximal and pure of width $k$.
The Tate twist $\W(1)$ is a variation of Hodge structures
of weight $k-2$. If we denote by $(E_1,\theta_1)$ the corresponding Higgs field,
and by $h_1^{p-1,q-1}$ the rank of $E_1^{p-1,q-1}$, then
\begin{gather*}
E^{k-\nu-1,\nu+1}=E_1^{k-2-\nu,\nu}\oplus \theta^{(\nu+1)}(E^{k,0})\otimes
\Omega^1_Y(\log S)^{-\nu-1}\\
\text{and \ \ }
h_1^{k-2-\nu,\nu}= h^{k-\nu-1,\nu+1}-h^{k,0},
\end{gather*}
hence
\begin{gather*}
\sum_{\nu=0}^{[\frac{k-2}{2}]} \deg(E_1^{k-2-\nu,\nu})=\hspace*{6cm}\\
\sum_{\nu=0}^{[\frac{k-2}{2}]} \big( \deg(E^{k-\nu-1,\nu+1}) -
(\nu+1) \cdot h^{k,0} \cdot \deg(\Omega_Y^1(\log
S))\big)=\\
\sum_{\nu=1}^{[\frac{k}{2}]} \big( \deg(E^{k-\nu,\nu})
-\nu \cdot h^{k,0} \cdot \deg(\Omega_Y^1(\log S))\big)=\\
\sum_{\nu=0}^{[\frac{k}{2}]} \big( \deg(E^{k-\nu,\nu})
-\nu \cdot h^{k,0} \cdot \deg(\Omega_Y^1(\log S))\big) -\deg(E^{k,0}).
\end{gather*}
Using d), (i), for $\nu=0$, and d), (ii), one obtains
\begin{gather*}
\sum_{\nu=0}^{[\frac{k-2}{2}]} \deg(E_1^{k-2-\nu,\nu})=\hspace*{6cm}\\
\sum_{\nu=0}^{[\frac{k}{2}]}
\frac{k-2\nu}{2}\cdot (h^{k-\nu,\nu}-h^{k,0}) \cdot \deg(\Omega^1_Y(\log S))
-\deg(E^{k,0})=\\
\sum_{\nu=1}^{[\frac{k}{2}]}
\frac{k-2\nu}{2}\cdot (h^{k-\nu,\nu}-h^{k,0}) \cdot \deg(\Omega^1_Y(\log S))=\\
\sum_{\nu=0}^{[\frac{k-2}{2}]}
\frac{k-2-2\nu}{2}\cdot h_1^{k-2-\nu,\nu} \cdot \deg(\Omega^1_Y(\log
S))=\\
\sum_{\nu=0}^{[\frac{k-2}{2}]}
\frac{k-2-2\nu}{2}\cdot h_1^{k-2-\nu,\nu}\cdot \deg(\Omega^1_Y(\log
S)).
\end{gather*}
Thereby $\W(1)$ satisfies again the equality d), (ii), hence
d), (i) as well, and by induction on the weight of the variation of Hodge structures one
obtains the decomposition asked for in \ref{gmaxhiggs} and \ref{maxhiggs}.
\end{proof}

\section{Decompositions of Higgs bundles with strictly maximal Higgs fields}
\label{secdecomp}
As in \cite{VZ4}, 1.2 and 1.4, one obtains as a corollary of
\cite{Si2}:

\begin{proposition}\label{semistable}
Assume that the Higgs field of a polarized $\C$ variation of Hodge
structures $\X$ of weight $k$ is strictly maximal. Then for each of the
strictly maximal pure components $F_i$ of width $i$
the sheaves $F_i^{p,k-p}$ are poly-stable.
\end{proposition}
\begin{proof}
Recall that $F_i^{p,k-p}\neq 0$ if and only if $k-i$ is even and
$k-i \leq 2p \leq k + i$.

Since the Higgs field is given by isomorphisms, it is sufficient
to show that for $\bar p = \frac{k+i}{2}$ the sheaf $F_i^{\bar
p,k-\bar p}$ is poly-stable. Let $\sA=\sA^{\bar p,k- \bar p}$ be a
sub-bundle and define inductively
$$
\sA^{p-1,k-p+1}=\tau_i (\sA^{p,k-p})\otimes
(\Omega^1_Y(\log S))^{-1} \subset F_i^{p-1,k-p+1}.
$$
Restricting the Higgs field, one finds that
$$
\big(\bigoplus_{p=\bar p -i}^{\bar p}\sA^{p,k-p}, \tau_i \big)
$$
is a sheaf, underlying a sub Higgs bundle of $(F_i,\tau_i)$. Hence
by \cite{Si2}
\begin{equation}\label{eq1}
\sum_{p=\bar p -i }^{\bar p}\deg(\sA^{p,k-p}) \leq
\sum_{p=\bar p -i}^{\bar p}\deg(F_i^{p,k-p})=0.
\end{equation}
On the other hand,
$$
\sum_{p=\bar p -i }^{\bar p}\deg(\sA^{p,k-p})=
i\cdot\deg(\sA)+ i\cdot \sum_{p=\bar p -i }^{\bar p}
\rank(\sA) \cdot \deg(\Omega^1_Y(\log S).
$$
Since the same equation holds true for
$F_i^{\bar p,k- \bar p}$ instead of $\sA$,
(\ref{eq1}) implies that
\begin{equation}\label{eq3}
\frac{\deg(\sA)}{\rank(\sA)} \leq
\frac{\deg(F_i^{\bar p,k- \bar p})}
{\rank(F_i^{\bar p,k- \bar p})}.
\end{equation}
Hence $F_i^{\bar p,k- \bar p}$ is semi-stable. If (\ref{eq3}) is
an equality,
$$
\big(\bigoplus_{p=\bar p -i}^{\bar p}\sA^{p,k-p}, \tau\big)
$$
must be a direct factor of $(F_i,\tau_i)$, hence $\sA$ is a direct
factor of $F_i^{\bar p,k- \bar p}$.
\end{proof}
\begin{remark}\label{pione}
Assume that $Y=\BP^1$ and that the Higgs field $(E,\theta)$ is non-zero.
\begin{enumerate}
\item[a.] The strict maximality of the Higgs field implies by \ref{semistable}
that
\begin{equation}\label{pione3}
F_i^{\frac{k+i}{2},\frac{k-i}{2}}=\bigoplus \sO_{\BP^1}(\nu)
\end{equation}
for some $\nu$, and that $2\cdot\nu=i\cdot (-2 + \#S)$.
If $k$ is odd, hence $F_i\neq 0$ for some odd number $i$, one finds
$\#S$ to be even.
\item[b.] Independently of $k$ one can always arrange the decomposition in
(\ref{pione3}) such that
$(F_i,\tau_i)=\bigoplus (N_i,\rho_i)$
where $N_i$ is a rank $i+1$ Higgs bundle with $N_i^{p,q}$ at most one
dimensional, and with a strictly maximal pure Higgs field.
\end{enumerate}
\end{remark}
\begin{remark}\label{pione2}
From now on we will frequently assume that $\#S$ is even. If $Y$ is a curve of genus $g>0$,
this holds true after replacing $Y$ by an \'etale covering of degree $2$.
For $Y=\BP^1$, and for $k$ odd, this assumption follows from the
strict maximality of the Higgs field. For $Y=\BP^1$ and $k$ even,
one may have to replace $Y$ by a twofold covering, ramified in
two points of $S$.
\end{remark}

If $\# S$ is even, we call an invertible sheaf $\sL$
a logarithmic theta characteristic, if there is an isomorphism
$$
\sL^2\> \simeq >> \Omega^1_Y(\log S) \mbox{ \ \ or equivalently \ \ }
\sigma:\sL \> \simeq >> \sL^{-1} \otimes \Omega^1_Y(\log S).
$$
So $(\sL\oplus \sL^{-1},\sigma)$ is an indecomposable
Higgs field of degree zero, hence it corresponds to a local $\C$
system $\LL$. This system is uniquely determined by $Y$, up to the
tensor product with local systems given by two division points in ${\rm Pic}(Y)$.

Let us fix once for all one of those local systems $\LL$. We will regard
$\LL$ as a $\C$ variation of Hodge structures of weight $1$ and width $1$, hence
concentrated in bidegrees $(1,0)$ and $(0,1)$. Then $S^i(\LL)$ has
weight $i$ and width $i$. Given $k$ with $k-i$ even, the Tate twist $S^i(\LL)
(-\frac{k-i}{2}) $ is a $\C$ variation of Hodge structures
of weight $k$ and width $i$.

The components of the Higgs bundle of
$S^i(\LL)(-\frac{k-i}{2})$ are concentrated in bidegrees $(p,k-p)$ for
$$
p=\frac{k-i}{2},  \frac{k-i}{2}+1,  \ldots , \frac{k+i}{2},
$$
and for $p=\frac{k+i-2\mu}{2}$ the corresponding component is isomorphic to
$\sL^{i-2\mu}$. The Higgs field is induced by $\sigma$, hence an isomorphism,
and we will denote it by
$$
S^{i-2\mu}(\sigma):\sL^{i-2\mu}\>>>
\sL^{i-2\mu-2}\otimes \Omega_Y^1(\log S).
$$
In particular $S^i(\LL) (-\frac{k-i}{2})$ is strictly maximal and pure.
\begin{proposition}\label{thetachar}
Assume again that $\X$ is a $\C$ variation of Hodge structures of
weight $k$ and with a strictly maximal Higgs field.
Assume in addition that $\#S$ is even.
Using the notations from \ref{maxpure}, let $\V_i$ be the strictly
maximal pure component of $\X$ of width $i$.
Then one has decompositions
\begin{equation}\label{dec4}
\V_i=S^i(\LL)(-\frac{k-i}{2}) \otimes \T_i,
\end{equation}
where $\T_i$ are unitary local system and a $\C$ variation of Hodge
structures concentrated in bidegree $(0,0)$.
\end{proposition}
\begin{proof} By assumption $\deg(\Omega_Y^1(\log S))$ is even,
hence that we can choose a logarithmic theta characteristic $\sL$.

Let $(F_i,\tau)$ denote the Higgs bundle corresponding to $\V_i$.
The strict maximality of the Higgs field
implies that for $\frac{k-i}{2} < p \leq \bar{p}=\frac{k + i}{2}$
$$
\tau^{p,k-p}:F_i^{p,k-p} \>>> F_i^{p-1,k-p+1}\otimes \Omega^1_Y(\log S)
$$
is an isomorphism. Let us write $g$ for the rank of those sheaves.
Then
$\det(F_i^{p,k-p})\simeq \det(F_i^{p-1,k-p+1})\otimes \sL^{2g}$,
and
$$
0=\sum_{p=\bar p -i}^{\bar p} \deg(F_i^{p,k-p})=
\sum_{\nu=0}^{i} \deg(F_i^{\bar p-i,k+i-\bar p})\cdot
2\nu g \deg(\sL).
$$
One obtains
\begin{equation}\label{eq4}
\deg(F_i^{\bar p -i ,k +i -\bar p}) = -i g \deg(\sL), \mbox{ \ \ and \ \ }
\deg(F_i^{\bar p,k-\bar p}) = i g \deg(\sL), .
\end{equation}
Choose $\sT=F_i^{\bar p,k-\bar p}\otimes \sL^{-i}$ with trivial
Higgs field. \ref{semistable} implies that this sheaf is
poly-stable, and by (\ref{eq4}) its degree is zero. Hence $(\sT,0)$
corresponds to a unitary $\C$ local system $\T$, which we may
regard as a $\C$ variation of polarized Hodge structures, concentrated in
bidegree $(0,0)$.

The $\bar p -\mu, k-\bar p + \mu$ component of the
Higgs bundle corresponding to the local system
$S^i(\LL)(- \frac{k-i}{2})\otimes \T$ is given by
$\sL^{i-2\mu} \otimes \sT$ and the Higgs field by
$S^{i-2\mu}(\sigma)\otimes {\rm id}_{\sT}$.

The isomorphism $\tau^{\bar p - \mu,k-\bar p + \mu}$, together with
$\Omega^1_Y(\log S)\simeq \sL^2$, induces isomorphisms
\begin{gather*}
\sT=F_i^{\bar p,k-\bar p}\otimes \sL^{-i} \> \delta_\mu >
\simeq >
F_i^{\bar p -\mu,k- \bar p +\mu}\otimes
\sL^{-i+ 2\mu)} \>\simeq >>\\
F_i^{\bar p -\mu -1,k- \bar p +\mu +1}\otimes \sL^{-i+ 2(\mu-1)}\>
\simeq >> F_i^{\bar p -i ,k+i -\bar p}\otimes \sL^{i}.
\end{gather*}
By the choice of the different isomorphisms
$$
\begin{CD}
\sL^{i-2\mu} \otimes \sT \> S^{i-2\mu}(\sigma)\otimes {\rm id}_\sT >>
\sL^{i-2(\mu-1)} \otimes \Omega^1_Y(\log S) \otimes \sT\\
\V {\rm id}_{\sL^{i-2\mu}}\otimes \delta_{\mu} VV
\V V {\rm id}_{\sL^{i-2(\mu-1)}\otimes \Omega_Y^1(\log S)}\otimes
\delta_{\mu-1} V\\
F_i^{\bar p -\mu, k- \bar p +\mu} \>>> \Omega_Y^1(\log S) \otimes
F_i^{\bar p -\mu -1 , k-\bar p + \mu + 1}
\end{CD}
$$
commutes, hence the two Higgs bundles are isomorphic.
\end{proof}
\begin{lemma}\label{tensorsplit}
Assume that $\X_{\bar\Q}$ is a variation of $\bar \Q\cap \R$ Hodge structures,
and that $\X=\X_{\bar\Q}\otimes \C$ has a strictly maximal Higgs field.
Assume moreover, that $\#S$ is even. Then, after replacing $Y$ by a finite
\'etale covering, the tensor product decomposition
in \ref{thetachar} can be defined over $\bar\Q$. To be more precise,
there exist a $\bar\Q\cap \R$ local system $\LL_{\bar\Q\cap \R}$,
and a unitary $\bar\Q\cap \R$ local system $\T_{i\bar\Q\cap \R}$ with
$$
\V_{i,\bar\Q} = S^i(\LL_{\bar\Q})(-\frac{k-i}{2})\otimes \T_{i\bar\Q}
\mbox{ \ \ and \ \ } \X_{\bar\Q} = \bigoplus_{i=0}^\ell \V_{i,\bar\Q}.
$$
\end{lemma}
\begin{proof}
By \ref{maxhiggs2} $\X_{\bar\Q} = \bigoplus_{i=0}^\ell \V_{i,\bar\Q}$
for some $\bar\Q$ variations of Hodge structures $\V_{i,\bar\Q}$
of width $i$, each one with a strictly maximal and pure Higgs field.

Hence we may restrict ourselves to the case that
$\X_{\bar\Q}=\V_{i,\bar\Q}$. Then the proof is similar to the one
given \cite{VZ4}, 3.7, ii), a) and b):\\

Consider for $j=\frac{k-i}{2}$
the isomorphism of local systems
$$\phi: S^i(\LL)(-j)\otimes \T_i \>
\simeq >> \V_i
$$
and the induced isomorphism
\begin{gather*}
\END_0(S^i(\LL)(-j) \otimes \T_i) = \hspace*{9cm}\\
\END_0(S^i(\LL)(-j)) \oplus \END_0(\T_i)\otimes \END_0(S^i(\LL)(-j))
\oplus \END_0(\T_i) \>\phi^2 >> \END(\V_i).
\end{gather*}
Since $\phi^2 \END_0(\T_i)$ is the unitary part of this
decomposition, by \cite{VZ4}, 3.3,
it is defined over $\bar{\Q}\cap\R$, as well as
$$
\phi^2(\END_0(S^i(\LL)(-j)) \oplus \END_0(\T_i)\otimes \END_0(S^i(
\LL)(-j)).
$$
The $j,-j$ part of the Higgs field corresponding to
$\phi^2\End_0(S^i(\LL)(-j))$ has rank one,
and its Higgs field is strictly maximal. Hence $\phi^2\End_0(S^i(\LL)(-j))$
is irreducible, and by \cite{VZ4}, 3.2, it is isomorphic to a
local system, defined over $\bar{\Q}$.
So
$$
\T_i\otimes \T_i \simeq \END(\T_i)\mbox{ \ \ and \ \ }
S^i(\LL)(-j) \otimes S^i(\LL)(-j)\simeq\END(S^i(\LL)(-j))
$$
are both isomorphic to local systems defined over some real
number field $K'$.

Consider for $\nu=2$ and $\mu=i$ or $\nu=g_0$ and $\mu=1$ the
moduli space $\sM (U, \Sl(\nu^{2\mu}))$ of reductive
representations of $\pi(U,*)$ into $\Sl(\nu^{2\mu})$. It is a
quasi-projective variety defined over $\Q$.
$$
S^i(\LL)(-j)\otimes S^i(\LL)(-j) \text{ \ \ (or} \T_i\otimes\T_i\text{)}
$$
is defined over $\bar\Q$, which implies that its isomorphism class in
$\sM(U, \Sl(\nu^{2\mu}))$ is a $\bar\Q$ valued point.

The morphism induced by the second tensor product
of the $\mu$-th symmetric product
$$
\rho:\sM(U, \Sl(\nu)) \>>> \sM(U, \Sl(\nu^{2\mu}))
$$
is clearly defined over $\Q$. As Simpson has shown (see
\cite{VZ4}, 3.4) $\rho$ is finite, hence the fibre
$$
\rho^{-1}([S^i(\LL)(-j)\otimes S^i(\LL)(-j)]) \text{ \ \ (or}
\rho^{-1}([\T_i\otimes \T_i])\text{)}
$$
consists of finitely many $\bar{\Q}$-valued points, hence $\LL$ and
$\T_i$ can both be defined over a number field $\bar\Q$.

Finally, passing to an \'etale
covering $\LL$ and $\T_i$ can both be assumed to be defined over
$\R$. In fact, the local system $\bar{\LL}$ has a strictly maximal Higgs
field, hence its Higgs field is of the form $(\sL'\oplus
{\sL'}^{-1}, \tau')$ where $\sL'$ is a theta characteristic.
Hence it differs from $\sL$ at most by the tensor product with a
two torsion point in ${\rm Pic}^0(Y)$. Replacing $Y$ by an
\'etale covering, we may assume $\LL=\bar{\LL}$.
Finally, if $(\sT_i,0)$ denotes the Higgs field of $\T_i$, and
$(\sT'_i,0)$ the one of $\bar\T_i$,
then $\sT_i\otimes \sL^i \simeq F_i^{j,k-j}$, whereas
$\sT'_i\otimes\sL^{-i}\simeq F_i^{k-j,j}$. Since the complex
conjugation interchanges $F_i^{j,k-j}$ and $F_i^{k-j,j}$
as well as $\sL$ and $\sL^{-1}$, one obtains that
$\bar{\T}=\T$.
\end{proof}

\section{Rigidity and the special Mumford Tate group}
\label{sectrigi}
Let us shortly recall the definition of the special
Mumford-Tate group (called Hodge group in \cite{Mum1} and \cite{Mum2},
see also \cite{De2} and \cite{Scho}). For a general fibre $F$ consider
the Hodge structure $H^k(F,\Q)$. The special
Mumford-Tate group $\Hg(F)$ is defined in as the
largest $\Q$ algebraic subgroup of $\Gl(H^k(F,\Q))$,
which leaves all Hodge tensors invariant,
i.e. all elements
$$
\eta\in \big[ \big(\bigotimes^{m}H^k(F,\Q)\big) \otimes
\big(\bigotimes^{m'}H^k(F,\Q)^\vee\big) \big]^{\frac{k(m-m')}{2},
\frac{k(m-m')}{2}}.
$$
$Hg(F)$ is a reductive group.

For a smooth family of manifolds $f:X_0 \to U$ with
$F=f^{-1}(y)$ for some $y\in U$, and for the corresponding $\Q$
variation of polarized Hodge structures $R^kf_*\Q_{X_0}$, consider
Hodge tensors $\eta$ on $F$ which remain
Hodge tensors under parallel transform. One defines the special
Mumford-Tate group
$\Hg(R^kf_*\Q_{X_0})$ as the largest $\Q$ subgroup which leaves
all those Hodge cycles invariant (\cite{De2}, \S 7, or
\cite{Scho}, 2.2).
\begin{lemma}\label{hg} \
\begin{enumerate}
\item[a.] $\Hg(R^kf_*\Q_{X_0})$ coincides for all $y'$ in an open
dense subset of $U$ with the special Mumford-Tate group
$\Hg(f^{-1}(y'))$. In particular it is reductive.
\item[b.] Let $G^{\rm Mon}$ denote the smallest reductive $\Q$ subgroup
of $\Gl(H^k(F,\R))$, containing the image $\Gamma$ of the monodromy
representation
$$
\gamma:\pi_0(U)\>>> \Gl(H^k(F,\R)).
$$
Then the connected component $G^{\rm Mon}_0$ of one inside
$G^{\rm Mon}$ is a subgroup of $\Hg(R^1f_*\Q_{X_0})$.
\item[c.] If $\#S$ is even and if $R^kf_*\C_{X_0}$ has a strictly maximal Higgs field
$$G^{\rm Mon}_0=\Hg(R^kf_*\Q_{X_0}).$$
\end{enumerate}
\end{lemma}
Recall that the assumption in \ref{hg}, c), implies that
the unitary part $\U=\V_0$ of $R^kf_*\C_{X_0}$
is either trivial, or for $k$ even concentrated in bidegree
$(\frac{k}{2},\frac{k}{2})$.
\begin{proof} a) has been verified in \cite{Scho}, 2.3.
As explained in \cite{De2}, \S 7, or \cite{Scho}, 2.4, the
Mumford-Tate group contains a subgroup of $\Gamma$ of finite index.
The same argument works for the special Mumford-Tate and b) holds true.

Since the special Mumford-Tate group is
reductive, part a) implies that $\Hg(R^kf_*\Q_{X_0})$ is reductive.
So
$$
G^{\rm Mon}_0 \subset \Hg(R^kf_*\Q_{X_0})
$$
is an inclusion of reductive groups. The proof of 3.1, (c), in
\cite{De3} carries over to show that both groups are equal,
if they leave the same tensors
$$
\eta\in \big(\bigotimes^{m}H^k(F,\Q)\big) \otimes
\big(\bigotimes^{m'}H^k(F,\Q)^\vee\big)
$$
invariant. Hence for c) one just has to verify, that such a tensor
$\eta$ can only be invariant under $G^{\rm Mon}_0$ or under $\Gamma$,
if it is of bidegree  $(\frac{k(m-m')}{2},\frac{k(m-m')}{2})$.

For $\eta\in \big(\bigotimes^{m}H^k(F,\Q)\big) \otimes
\big(\bigotimes^{m'}H^k(F,\Q)^\vee\big)$, invariant under
$\Gamma$, let $\tilde\eta$ be the corresponding global section of
$$
\V(m,m')_\Q= \big(\bigotimes^{m}(R^kf_*\Q_{X_0})\big) \otimes
\big(\bigotimes^{m'}(R^kf_*\Q_{X_0})^\vee\big).
$$
By \ref{thetachar} one has a decomposition
$$
R^kf_*\C_{X_0}=\bigoplus_{i = 0}^{\ell} \V_i \simeq
\bigoplus_{i=0}^\ell \big( S^i(\LL)(-\frac{k-i}{2}) \otimes
\T_i \big),
$$
where $\T_i$ is a unitary local system and a $\C$ variation of Hodge
structures concentrated in bidegree $(0,0)$, and where
$S^i(\LL)(-\frac{k-i}{2})$ has a
strictly maximal Higgs field of width $i$.
In order to prove c) it remains to verify for
$$\V(m,m')=\V(m,m')_\Q\otimes_\Q\C:$$
\begin{claim}\label{ppvan}
$H_0(U,\V(m,m'))^{p,q}=0$ for $(p,q)\neq (\frac{k(m-m')}{2},
\frac{k(m-m')}{2})$.
\end{claim}
\noindent
{\it Proof.} \
By \cite{FH}, page 80, $\V(m,m')$ decomposes in a direct sum of local
sub systems, each of which is a tensor product of Schur functors
$\Shur_\lambda$, for certain partitions $\lambda$, applied to
one of the local systems $\U$, $\U^\vee$, $\LL$, $\LL^\vee$,
$\T_i$, and $\T_i^\vee$, for $i= 1, \ldots , \ell$.

Since $\LL$ is of rank two, and since $\LL$ is self dual,
the only possible Schur functors with $\Shur_\lambda(\LL)$
or $\Shur_\lambda(\LL^\vee)$ non zero, are isomorphic to $S^\nu(\LL)$,
for $\nu \geq 0$. They all have a strictly maximal and pure Higgs field of
width $>0$, except of $S^0(\LL)\simeq (\wedge^2(\LL))^a$. In fact,
as a variation of Hodge structures, the latter will
always be concentrated in some bidegree $(\iota_1,\iota_1)$.

Note furthermore that all $\Shur_{\lambda}(\U)$,
$\Shur_{\lambda}(\U^\vee)$, $\Shur_{\lambda}(\T_i)$ or
$\Shur_\lambda(\T^\vee)$
are unitary of bidegree $(\iota_2,\iota_2)$, for some $\iota_2$.
The same remains true for all their tensor products.

Consider a typical direct factor
$$
\M=S^{\nu}(\LL) \otimes \Shur_{\lambda_0}(\U)
\otimes \Shur_{\lambda'_0} \otimes \bigotimes_{i=1}^{\ell} \big(
\Shur_{\lambda_i}(\T_i)\otimes \Shur_{\lambda'_i}(\T^\vee_i)\big)
$$
of $\V(m,m')$. Taking the tensor product of a local system with
a strictly maximal pure Higgs field with a unitary local system,
one obtains again a local system with a strictly maximal pure Higgs field.
Since a variation of Hodge structures of weight $>0$ and
with a strictly maximal pure Higgs field can not have
any global section, one finds that $H^0(U,\M)=0$, except for $\nu=0$.
In the latter case, seen as a variation of
Hodge structures, $\M$ is concentrated in bidegree $(\iota,\iota)$, hence
$H^0(U,\M)$, as well.
\end{proof}
\begin{lemma}\label{rigid} Let $f:X\to Y$ be a non-isotrivial
semi-stable family of $n$-dimensional varieties, with $X_0\to U=Y\setminus S$
smooth. Assume that for a general fibre $F$ of $f$ the local Torelli
theorem holds true. If $R^nf_*\C_{X_0}$ has a strictly maximal
Higgs field, then $f:X \to Y$ is rigid, i.e. there exists no non-trivial
smooth deformations $\tilde f : \tilde X_0 \to U\times T$ of $f:X_0\to U$.
\end{lemma}
\begin{proof} By \cite{Fal} (see also \cite{Pet}) the rigidity follows from
the vanishing of
$$
\End(R^nf_*\C_{X_0})^{-1,1}.
$$
Using the notation introduced in the proof of \ref{hg}, this group is
$$
H^0(U,\V(1,1))^{-1,1} \text{ \ \ for \ \ }
\V(1,1)= R^nf_*\C_{X_0}\otimes R^nf_*\C_{X_0}^{\vee},
$$
hence zero by \ref{ppvan}.
\end{proof}

\section{Splitting over $\Q$ for $S\neq \emptyset$}\label{secsplit}
As in \cite{VZ4} we will show in this section that, replacing $Y$ by an
\'etale covering (or by a degree two covering ramified in two
points for $Y=\BP^1$), the local system $\LL$ can be defined over
$\Q$ and that the unitary local systems $\T_i$ in \ref{thetachar}
can be assumed to be trivial, i.e. a direct sum of copies of $\C$.
We will need the following well known property of the trace.
\begin{proposition}\label{trace}
Fixing a natural number $n\in\N$, there exists
a polynomial $p_n(t)\in \Z[t]$ of degree $n$ and with the leading
coefficient $1$, such that $${\rm tr}({\rm Sym}^n(M))=p_n({\rm tr}(M)),
\mbox{ \ \ \ for all \ \ \ } M\in  \Sl_2(\C).$$
\end{proposition}
\begin{proof} We may assume $M$ is of Jordan normal form
$$\left [\begin{array}{cccc}
\lambda & \alpha\\
0        &\lambda^{-1}\\
\end{array}\right],$$
which acts on a two dimensional vector space $V$ over $\mathbb C$ with a
basis $\{x, y\}$ by
$$\left [\begin{array}{cccc}
x\\
y\\
\end{array}\right]\mapsto  \left [\begin{array}{cccc}
\lambda & \alpha\\
0        &\lambda^{-1}\\
\end{array}\right] \left [\begin{array}{cccc}
x\\
y\\
\end{array}\right].$$
Using the basis $\{x^n,\, x^{n-1}y,\, x^{n-2}y^2,..., y^n\}$ for
${\rm Sym}^n V$ one finds $$ {\rm tr}({\rm
Sym}^n(M))=\lambda^n+\lambda^{n-2}+\lambda^{n-4}+...+\lambda^{-n}.$$
It is easy to see there exist integers $a_{n-1},\,
a_{n-2},\,...,a_0\in\mathbb Z$ such that $\forall
\lambda\in\mathbb C^*$
\begin{gather*}
\lambda^n+\lambda^{n-2}+\lambda^{n-4}+...+\lambda^{-n}=\\
(\lambda+\lambda^{-1}
) ^n +
a_{n-1}(\lambda^{n-1}+\lambda^{(n-1)-2}+...+\lambda^{-(n-1)})+\\
a_{n-2}(\lambda^{n-2}+\lambda^{(n-2)-2}+...+\lambda^{-(n-2)})+ \cdots
+a_0=\\
({\rm tr}(M))^n+a_{n-1}{\rm tr}({\rm Sym}^{n-1}(M))+a_{n-2}{\rm tr}({\rm
Sym}^{n-2}(M))+...+a_0. \end{gather*}
The proposition follows by induction on $n.$
\end{proof}
\begin{theorem}\label{nonempty}
Assume that $S=Y\setminus U\neq \emptyset$, that $\#S$ is even, and that
$\X_\Q$ is a non-constant polarized $\Q$ variation of Hodge structures of weight
$k$ with unipotent monodromies around all $s\in S$. If $\X=\X_\Q\otimes\C$
has a strictly maximal Higgs field, then, after replacing $Y$
by an \'etale covering,
there exists a rank $2$ polarized $\Z-$ variation of Hodge structures $\LL$
with an $\Q$-Hodge isometry
\begin{equation}\label{decnonempty}
\X_\Q \simeq \bigoplus_{i=0}^{[\frac{k}{2}]} S^{k-2i}(\LL_\Q)(-i)^{\oplus h^{k-i,i}-
h^{k-i+1,i-1}},
\end{equation}
where $h^{k-i,i}$ denote the Hodge numbers of $\X$.
\end{theorem}
Recall that the definition of a strictly maximal Higgs field implies that
$$h^{k-i+1,i+1} \leq h^{k-i,i},$$
for $i=0,\ldots,[\frac{k}{2}]$. So the exponents in \ref{decnonempty} are all
non negative.
\begin{proof}[Proof of \ref{nonempty}]
By Definition \ref{gmaxhiggs} and by Proposition \ref{thetachar}
we have a decomposition
\begin{equation}\label{decomp1}
\X=\X_\Q\otimes\C = \bigoplus_{i=0}^{[\frac{k}{2}]}\V_{k-2i}=
\bigoplus_{i=0}^{[\frac{k}{2}]}S^{k-2i}(\LL)(-i) \otimes \T_{k-2i}.
\end{equation}
The local system $\X\otimes\X$ decomposes in a direct sum
\begin{equation}\label{decomp2}
\V_{k-2i}\otimes \V_{k-2j}= S^{k-2i}(\LL)(-i)\otimes S^{k-2j}(\LL)(-j)\otimes \T_{k-2i}
\otimes \T_{k-2j}
\end{equation}
for all pair $i,j\in\{0,\ldots,[\frac{k}{2}]\}$. Since $\det(\LL)=\C$, for $i \leq j$
Pieri's formula (see \cite{FH}, 6.16, and 6.9) implies that
\begin{gather}\label{decomp3}
S^{k-2i}(\LL)(-i)\otimes S^{k-2j}(\LL)(-j)=
S^{2k-2i-2j}(\LL)(-i-j) \oplus\\ S^{2k-2i-2j-2}(\LL)(-i-j-1) \oplus \cdots \oplus
S^{-2i+2j}(\LL)(-i+j).\notag
\end{gather}
For $a>0$ the polarized variation of Hodge structures $S^{2a}(\LL)(-a)$ have a
strictly maximal pure Higgs field, as well as their tensor product with a
unitary local system, concentrated in one bidegree $(\iota,\iota)$.
So the only unitary sub local systems of $\X\otimes\X$ are contained in
$\V_{k-2i}\otimes\V_{k-2i}$ and they are of the form
$\T_{k-2i}\otimes \T_{k-2i}$. By \cite{VZ4}, 3.3,
$$
\T=\bigoplus_{i=0}^{[\frac{k}{2}]}\T_{k-2i}\otimes \T_{k-2i}
$$
comes from a direct factor of $\X_K=\X_\Q\otimes K$ for some number field $K$.
Since  $S\not=\emptyset,$ the argument used to prove \cite {VZ4}, 4.1, carries
over to show that on can even choose $K=\Q$, and that $\T$ admits an $\Z$-structure.
Since  the eigenvalues of the local monodromies of $ \X$ around $S$ are all equal to 1,
the local monodromies of $\T$ around $S$ are all identity. Hence,
$\T$ extends to a unitary local system on $Y$, say of rank $\kappa$.
The representation defining $\T$ can be written as
$$
\rho:\pi_1(U,*) \to \pi_1(Y,*) \to \Gl(\kappa,\Q)
$$
and $\rho$ factors through a finite quotient of $\pi_1(Y,*)$. Replacing
$Y$ by a finite covering, we may assume that $\rho$ is trivial, hence that
$\T_{k-2i}\otimes \T_{k-2i}$ is trivial for all $i$.
Writing $\kappa_{k-2i}$ for the rank of $\T_{k-2i}$, hence
$$
\kappa_{k-2i}=h^{k-i,i}-h^{k-i+1,i-1}
$$
the tensor product induces a map
$$
\Gl(\kappa_{k-2i},\C)\to \Gl(\kappa^2_{k-2i},\C)
$$
with finite kernel, hence we may assume, replacing $Y$ again by some \'etale
cover, that $\T_{k-2i}$ is trivial, for all $i$.\\

Let us consider again the decompositions in (\ref{decomp1}), (\ref{decomp2}) and
\ref{decomp3}). They imply, that for some $n>0$
the local system $S^n(\LL)$ is a sub system of $\X\times \X$ which is
defined over $\Q$.

We need to show $\LL$ is isomorphic to a rank two local system admitting an
$\Z$-structure. Using the trace formula in Proposition \ref{trace} one sees that
${\rm tr}(\LL)$ is contained the ring of algebraic integers $\sO_K$ of an
algebraic number field $K\subset \C$. By  the argument used in the proof of Proposition
5.5 in \cite{VZ4} one shows further that $\sigma ({\rm tr}(\LL))$ is bounded for
any embedding $\sigma: K\to\C$ unequal to the identity.
So applying Takeuchi's Theorem
(see Theorem 5.3 and Proposition 5.5 in \cite {VZ4}) one obtains:
\begin{enumerate}
\item[i.] $\LL$ is  derived from a quaternion algebra $A$
over a totally real number field $F$ with $d$ distinct embeddings
$$\sigma_1={\rm id},\,\sigma_2,\ldots,\sigma_d: F\>>> \R.$$
\item[ii.] For $1\leq i\leq d$ there exists $\R$-isomorphism
$$\rho_1: A^{\sigma_1}\otimes\R\simeq M(2,\R),\mbox{ \ and \ }
\rho_i: A^{\sigma_i}\otimes\R\simeq \mathbb H,\mbox{ \ for \ }
2\leq i\leq d.$$
\item[iii.] The representation $\Gamma_{\LL}  :\pi_1(Y,*)\to \Sl(2,\R)$
defining the local system $\LL$ factors like
$$
\pi_1(Y,*)\>\simeq >> \Gamma \subset \rho_1(\sO^1) \>>>
\Sl(2,\R\cap\bar\Q)\subset \Sl(2,\R),
$$
and $Y\simeq \sH/\Gamma$, where $\sO^1$ is the group of units in an order
$\sO\subset A$ over $F$, and where $\sH$ denotes the upper half plane.\\
\item[iv.] $F(\sqrt{a})$ is a field of definition for $\LL$, for some $a\in
F$.
\item[v.] If $\tau_i,\,1\leq i\leq d$ are extensions of $\sigma_i$
to $F(\sqrt{a})$, and if $\LL_i$ denotes the local system defined
by
$$
\pi_1(Y,*)\>>> \Sl(2,F(\sqrt{a}))\> \tau_i>> \Sl(2,\bar\Q),
$$
then $\LL_i$ is a unitary local system, for $i>1$, and $\LL_1\simeq \LL$.
\item[vi.] Up to isomorphism, $\LL_i$ does not depend on
the extension $\tau_i$ chosen.
\end{enumerate}
Using the notations introduced above, the condition $S\neq \emptyset$ implies that
$\sH/\Gamma$ is not projective. This is only possible for $F=\Q$ and
for $A=M(2,\Q)$. Then, replacing $Y$ by a finite \'etale
cover, $\LL$ is defined over $\Z.$
\end{proof}

\section{Variations of odd weight}
\label{sectodd}
In this section we will show, that for $k$ odd, and $S=\emptyset$ there are no
polarized $\Q$ variations of Hodge structures $\X_\Q$ of weight $k$
with a strictly maximal Higgs field
and with Hodge number $h^{k,0}=1$. In fact, we will only use the highest width
part of $\X$ to find a contradiction.

\begin{proposition}\label{empty}
Over a projective curve $Y$ there exists no polarized $\Q$ variation
of Hodge structures of odd weight $k$ which allows a decomposition
$$
\X=\X_\Q\otimes\C=\V_k\oplus \W
$$
of polarized $\C$ variations of Hodge structures with:
\begin{enumerate}
\item[(i)] The highest Hodge number $h^{k,0}=1$.
\item[(ii)] The width of $\W$ is strictly smaller than $k$.
\item[(iii)] $\V_k$ is strictly maximal and pure
of width $k$.
\end{enumerate}
\end{proposition}
If $\X_\Q$ is a polarized $\Q$ variation of Hodge structures with a strictly maximal Higgs
field, and with highest Hodge number $h^{k,0}=1$, then the assumptions in \ref{empty}
hold true by Lemma \ref{deg}. So we obtain
\begin{corollary}\label{empty2}
Over a projective curve there exists no non-constant polarized $\Q$ variation of Hodge structures
of odd weight $k$, with a strictly maximal Higgs field and with $h^{k,0}=1$.
\end{corollary}
\begin{proof}[Proof of \ref{empty}]
Assume there exists an $\X_\Q$, with $k$ odd satisfying the conditions i), ii)
and iii).
By \ref{thetachar}, applied to $\V_k$ instead of $\X$, one has a decomposition
$$\V_k\simeq S^k(\LL)\otimes \T_k$$
with $\LL$ a rank $2$ weight $2$ variation of Hodge structures
with a strictly maximal pure Higgs field, and with $\T_k$ a unitary local system.
The assumption $h^{k,0}=1$ implies that $T_k$ is of rank one, hence
trivial.

Let $\sigma$ be an automorphism of $\bar\Q$. By \cite{Si2} (see also
\cite{VZ4}, 5.5) $\LL^\sigma$ is again a $\C$ variation of Hodge
structures of weight $2$. Hence it either is unitary, or it has a
strict generically maximal Higgs field.

If $\LL^\sigma$ is unitary $\V_k^\sigma \simeq S^k(\LL^\sigma)$ is unitary, hence
$\X$ must have a unitary part. This however, for $k$ odd, is excluded by the
strict maximality of the Higgs field.

Hence $\LL^\sigma$ has again a strict generically maximal Higgs field, and
the width of $\V_k^\sigma\simeq S^k(\LL^\sigma)$ is $k$. This implies
that the image of $\V_k^\sigma$ under the projection
$\X \to \W$ is zero. The same holds true for the image of $\W^\sigma$
under the projection $\X \to \V_k$. Hence $\V_k^\sigma = \V_k$ and
$\W^\sigma=\W$, as sub local systems of $\X$.
Then $\LL^\sigma$ has a strictly maximal Higgs field, and after passing to
some \'etale covering it is isomorphic to $\LL$.

Since $\V_k^\sigma = \V_k$, the sub local system $\V_k$ of $\X$ is defined over $\Q$,
Hence $S^k(\LL)$ as well. Recall that for $a \geq b$ the Pieri formula
(see \cite{FH}, 6.16, and 6.9)
$$
S^a(\LL)\otimes S^b(\LL)=S^{a+b}(\LL)\oplus S^{a+b-2}(\LL)....\oplus
S^{a-b}(\LL).
$$
Applying $\sigma$ and using again that $S^i(\LL)$ is strictly maximal and pure
of width $i$, one finds that $\sigma$ respects this
decomposition. In particular, choosing $a=b=k$ one finds $S^2(\LL)$ as a sub
local system of $S^k(\LL)\otimes S^k(\LL)$ invariant under all automorphisms
of $\bar\Q$. Hence
$$S^{2k-2[\frac{k}{2}]}(\LL)$$
is defined over $\Q$. Next one repeats this argument for
$a=2k-2[\frac{k}{2}]$ and $b=k$, and one finds for $a-b=k-2[\frac{k}{2}]=1$
that $\LL$ is defined over $\Q$.

So $\LL$ admits an $\Z-$structure and there exists some
$\Gamma\subset \Sl_2(\Z)$ with $Y\simeq \sH/\Gamma$.
Then however $Y$ can not be projective.
\end{proof}

\section{Families of Calabi-Yau manifolds and other examples}\label{secex}
C. Borcea constructed in \cite{Bor} Calabi-Yau threefolds
$F$ as a resolution of quotients
$$
E_1\times E_2 \times E_3/G
$$
where $E_i$ are elliptic curves with involution $\iota_i$ and
$$
G=\Z_2\times\Z_2 = <{\rm id}\times\iota_2\times\iota_3, \
\iota_1 \times{\rm id}\times\iota_3>.
$$
For $E_1=E_2=E_3$ this construction easily extends to modular families
of elliptic curves, and to higher dimensions.
Remark that the group $G$ is a subgroup of index two of the group generated by
the involutions $\iota_i$ acting on the $i$-th factor. Hence
$F$ can also be obtained as a twofold covering of some blowing up of
$(\BP^1)^3$.
\begin{example}\label{example3}
For all $n \geq 1$ there exists families $f:X\to Y$ of Calabi-Yau $n$-folds,
with smooth part $X_0\to U$, with $S=Y\setminus U\neq \emptyset$,
such that $R^nf_*\Q_{X_0}$ has a strictly maximal Higgs field.
\end{example}
\begin{proof}
Let $h: E \to Y$ be a modular family of elliptic curves, smooth over
$U\subset Y$ and birational to a semistable family. As explained in
\cite{VZ4}, Section 2, $R^1h_*\Q_{E_0}$ has a strictly maximal Higgs field,
where $E_0=h^{-1}(U)\to U$ is the smooth part.

We will assume, for simplicity, that the subgroup $E_{(2)}$
of two-division points in $E$ is trivial, i.e. that it consist
of the union of $4$ sections of $h$, pairwise disjoint on $E_0$.
Blowing up the singular fibres, we may assume
that the involution $\iota_0$ of $E_0\to U$ extends to an
involution $\iota$ on $E$, and that the quotient $\BP=E/\iota$ is non singular.
We will write $\varphi:\BP \to Y$ for the induced morphism
and $\Pi_{(2)}$ for the image of $E_{(2)}$ in $\BP$.

For any $n >0$ consider the $n$-fold product $E^n=E\times_Y \cdots \times_Y E$
and the group $H$ of automorphisms of $E^n$, generated by the involutions
$\iota$ on the different factors. $G$ again denotes the subgroup generated
by automorphisms $\eta_{ij}$ with $1\leq i < j \leq n$,
whose restriction to the $k$-th factor of $E^n$ is the identity, except for
the $i$-th or $j$-th factor, where it acts as $\iota$. So $E^n/G$
is a two to one covering of
$$
E^n/H=\BP^{(n)}=\BP \times_Y \cdots \times_Y\BP \> \varphi^{(n)} >> Y.
$$
The subgroup of $H$ leaving a components of $E_{(2)} \times_Y E^{n-1}$
pointwise invariant is generated by $\iota\times {\rm id}^{n-1}$, hence
$\Pi_{(2)}\times_Y \BP^{(n-1)}$ is totally ramified in $E^n/G$. Up to permuting
the factors, each component of the discriminant divisor $\Gamma \subset \BP^{(n)}$
is of this form. Let
$$\big\{ \Upsilon_i ; \ i\in \{1,\ldots,M \}\big\}
$$
be the set of $2$ by $2$ intersections of components of $\Gamma$.
To choose an ordering one may start, for example,
with the $16$ intersections of components given by two torsion points on
the first and second factor, then those coming from the first and third one,
and so on. Consider the sequence of blowing ups
$$
\BP^{(n,M)} \> \gamma_M >> \BP^{(n,M-1)} \>\gamma_{M-1} >> \cdots \>\gamma_2 >>
\BP^{(n,1)} \>\gamma_1 >> \BP^{(n,0)}=\BP^{(n)},
$$
where $\gamma_i$ the blowing up of the proper transform of $\Upsilon_i$ in
$\BP^{(n,i-1)}$. Write
$$
\gamma=\gamma_M \circ \cdots \circ \gamma_1:\BP^{(n,M)} = \tilde \BP^{(n)}
\>>> \BP^{(n)}
$$
for the composite. We choose $X$ to be the normalization of $\tilde \BP^{(n)}$
in the function field of $E^n/G$. We write $f:X\to Y$
and $\delta:X\to E^n/G$ for the induced morphisms.
\begin{claim}\label{cl1ex3}
The morphism
$$
f:X_0=f^{-1}(U) \to U
$$
is a smooth family of Calabi-Yau manifolds.
\end{claim}
\begin{proof}  Consider for $y\in U$ the fibre
$\pi_y:F_y=f^{-1}(y) \to \tilde \BP^{(n)}_y$.
The twofold covering $E^n_y/G \to \BP_y^{(n)}=(\BP^1)^n$ is given
by the invertible sheaf
$$
\sL_y=\sO_{(\BP^1)^n}(2,\ldots,2)
$$
and a section of $\sL_y^2$ with divisor $\Gamma_y$. The morphism $\gamma_y :\tilde
\BP_y^{(n)} \to \BP_y^{(n)}$ is the blowing up of all $2$ by $2$ intersections
of the normal crossing divisor $\Gamma_y$. So $\gamma_y^*\Gamma_y$ is again a
normal crossing divisor.

The proper transform $\tilde \Gamma_y$ of $\Gamma_y$ in
$\tilde \BP^{(n)}_y$ consists of disjoint irreducible components, and
$$
\gamma_y^*\Gamma_y-\tilde \Gamma_y=2\cdot \Phi_y
$$
for some reduced divisor $\Phi_y$. Since the covering
$\pi_y:F_y=f^{-1}(y) \to \tilde \BP^{(n)}_y$ is obtained
by taking the second root out of the divisor
$\gamma_y^*\Gamma_y$, the discriminant of $\pi_y$ is $\tilde \Gamma_y$.
The latter being non singular implies that $F_y$ is non singular.
Moreover, for all $p$, one has
\begin{equation}\label{gl2ex3}
\pi_{y*}\Omega^p_{F_y}= \Omega^p_{\tilde \BP^{(n)}_y} \oplus
\gamma_y^*\sL_y^{-1}\otimes\sO_{\tilde \BP^{(n)}_y}( \Phi_y )\otimes
\Omega^p_{\tilde \BP^{(n)}_y}(\log(\tilde \Gamma_y)).
\end{equation}
For $p>0$ the first direct factor can not have any global section.
The second one is a subsheaf of
$$
\gamma_y^*(\sL_y^{-1}\otimes \Omega^p_{\BP^{(n)}_y}
(\log \Gamma_y))\otimes\sO_{\tilde \BP^{(n)}_y}( \Phi_y ).
$$
Since $\Phi_y$ is contained in the exceptional locus of $\gamma_y$,
and since the ampleness of $\sL_y$ implies that for $p<n$
$$
H^0(\BP^{(n)}, \sL_y^{-1}\otimes \Omega^p_{\BP^{(n)}_y}(\log \Gamma_y))=0,
$$
one finds $H^0(F_y,\Omega^p_{F_y})=0$ for $0 < p < n$.
Finally, for $p=n$,
\begin{gather*}
\pi_{y*}\omega_{F_y}= \omega_{\tilde \BP^{(n)}_y} \oplus
\gamma_y^*\sL_y^{-1}(\Phi_y)\otimes
\omega_{\tilde \BP^{(n)}_y}(\log(\tilde \Gamma_y))=\\
\omega_{\tilde \BP^{(n)}_y} \oplus
\gamma_y^*(\sL_y^{1}\otimes \omega_{(\BP^1){n}})=\sO_{\tilde \BP^{(n)}_y},
\end{gather*}
and $\omega_{F_y}$ has a nowhere vanishing section.
\end{proof}
The description of the covering $X\to \tilde \BP^{(n)}$ given above carries over
to the whole family.

In the sequel ${}_0$ always denotes the preimage of $U$,
in particular $X_0=f^{-1}(U)$, $E_0=h^{-1}(U)$, and $\BP^{(n)}_0$
is a $(\BP^1)^n$ bundle over $U$. We choose
$\sL_0$ to be an invertible sheaf which fibrewise
coincides with $\sO_{(\BP^1)^n}(2,\ldots,2)$
and with $\sL_0^2=\sO_{\BP^{(n)}_0}(\Gamma_0)$.

It remains to show that $R^n f_* (\Q_{X_0})$ has a strictly maximal Higgs
field. A generator $\iota'$ of $H/G$ act on $R^n f_* (\Q_{X_0})$,
and we write $R^n f_* (\Q_{X_0})^{(1)}$ for the invariants
and $R^n f_* (\Q_{X_0})^{(-1)}$ for the anti-invariants.
Then the strict maximality of the Higgs field
follow from:
\begin{claim}\label{cl2ex3}
There is an isomorphism of variation of Hodge structures
$$
\bigotimes^n R^1 h_*(\Q_{E_0}) \>\theta>\simeq > R^n f_* (\Q_{X_0})^{(-1)}.
$$
Moreover, $R^n f_* (\Q_{X_0})^{(1)}=0$ for $n$ odd, and
$R^n f_* (\Q_{X_0})^{(1)}$ is unitary of bidegree $(\frac{n}{2},\frac{n}{2})$,
for $n$ even.
\end{claim}
\noindent{\it Proof.} \
The K\"unneth decomposition and the trace give natural maps
$$
\bigotimes^n R^1 h_*(\Q_{E_0}) \>>> R^n h^n_* (\Q_{E^n_0})
\>>> R^n f_* (\Q_{X_0}).
$$
For the invariants $R^n f_* (\Q_{X_0})^{(1)}$ the Higgs bundle is
\begin{gather*}
\bigoplus_{p+q=n} R^q \varphi^{(n)}_* \Omega^p_{\tilde \BP_0^{(n)}}.
\end{gather*}
Remark that the local systems for $\BP_0^{(n)} \to Y$ are trivial local systems,
and zero for $n$ odd. Assume the same for $\BP_0^{(n,i-1)} \to Y$. By construction,
$\BP_0^{(n,i)} \to Y$ is obtained by blowing up the proper transform
$\Upsilon'_i$ of $\Upsilon_i$ in $\BP_0^{(n,i-1)}$. The exceptional divisor is
a $\BP^1$ bundle over $\Upsilon'_i$, and the latter is itself a locally constant
family of rational varieties. Hence again the $n$-th direct
image of $\Q_{\BP_0^{(n,i)}}$ on  $U$ is a trivial local system,
and zero for $n$ odd. By induction one obtains the second part of \ref{cl2ex3}.

For the anti-invariants $R^n f_* (\Q_{X_0})^{(-1)}$ the Higgs bundle is the
direct sum of all
$$
F_0^{p,q}=R^q \varphi^{(n)}_*\gamma_*
\Omega^p_{\tilde \BP_0^{(n)}/U}(\log \tilde \Gamma_0)\otimes \gamma^* \sL_0^{-1}
\otimes \sO_{\tilde \BP_0^{(n)}}(\Phi_0),
$$
with $p+q=n$. The sheaf $F_0^{p,q}$ contains
$$
R^q \varphi^{(n)}_*\gamma_* \Omega^p_{\tilde \BP_0^{(n)}/U}(\log \gamma^*
\Gamma_0)\otimes \gamma^* \sL_0^{-1},
$$
and the cokernel is
\begin{equation}\label{gl3ex3}
R^q \varphi^{(n)}_*\gamma_*
\Omega^{p-1}_{\Phi_0/U}(\log \tilde \Gamma_0 \cap \Phi_0) \otimes \gamma^*
\sL_0^{-1} \otimes \sO_{\Phi_0}(\Phi_0).
\end{equation}
As above, each component of $\Phi_0$ is a $\BP^1$ bundle over some
$\Upsilon'_i$ and $\Gamma_0\cap \Phi_0$ consists of two disjoint sections
of this bundle. Since $\sO_{\Phi_0}(\Phi_0)$ restricted to its
fibres is $\sO_{\BP^1}(-1)$, the sheaf in (\ref{gl3ex3}) is zero,
for all $p$ and $q$.

This implies that
$$
F_0^{p,q}=R^q \varphi^{(n)}_*\Omega^p_{\BP_0^{(n)}/U}(\log \Gamma_0)\otimes \sL_0^{-1}.
$$
On the other hand,
$$
R^q \varphi_*\Omega^p_{\BP_0/U}(\log \Pi_{(2)})\otimes \sO_{\BP_0}(-2)=0
$$
is zero for $(p,q)\neq (0,1)$ and $(p,q)\neq (1,0)$, and
$$
\bigoplus_{p+q=n} F_0^{p,q}=
R^n \varphi^{(n)}_*\Omega^\bullet_{\BP_0^{(n)}/U}(\log \Gamma_0)\otimes \sL_0^{-1}
=\bigotimes^n R^1 \varphi_* \Omega^\bullet_{\BP_0/U}(\log \Pi_{(2)})
\otimes \sO_{\BP_0}(-2),
$$
as claimed.
\end{proof}
\begin{remark}
Of course the proof of \ref{cl1ex3} and \ref{cl2ex3} allows to determine
the invariants of the fibre $F$. Let us just remark that for $n=3$, i.e.
for Borcea's example, $f:X_0\to U$ is a family of Calabi Yau manifolds
with Euler number $e=96$ and $h^{2,1}=3$. One has
$$R^3f_*(\Q_{X_0})=\bigotimes^3R^1h_*(\Q_{E_0})=S^3R^1h_*(\Q_{E_0})\oplus
R^1h_*(\Q_{E_0})^{\oplus 2}.$$
\end{remark}
Starting from families of Abelian varieties with strictly maximal Higgs fields,
as classified in \cite{VZ4}, Theorem 0.5, it is easy to give
examples of smooth families $f:X\to Y$ of $n$-dimensional manifolds
over a projective curve $Y$, with a strictly maximal Higgs field.

Let us first recall the classification of smooth families of Abelian varieties
with a strictly maximal Higgs field:\\

Let $A$ be a quaternion division algebra defined over a totally real
number field $F$, which is ramified at all infinite places
except one. Choose an embedding
$$
D=\Cor_{F/\Q}A \subset M(2^m,\Q),
$$
with $m$ minimal, hence for $d=[F:\Q]$ one has $m=d$ or $m=d+1$.
Then the moduli functor of Abelian varieties with special Mumford-Tate group
$$
\Hg=\{ x\in D^*; \ x{\bar x}=1\}
$$
and with a suitable level structure, is represented by a
smooth family $Z_A \to Y_A $ over a compact Shimura curve $Y_A$.
As explained in \cite{VZ4} $Z_A \to Y_A$ has a strictly maximal Higgs field,
and the general fibre $Z_\eta$ of $Z_A \to Y_A$ is of the following type:
\begin{enumerate}
\item[i.] $1< m=d$ odd. In this case
$\dim(Z_\eta)=2^{d-1}$ and ${\rm End}(Z_\eta)\otimes_\Z \Q = \Q$.
\item[ii.] $m=d+1$. Then $\dim(Z_\eta)=2^{d}$ and
\begin{enumerate}
\item[a.] for $d$ odd, ${\rm End}(Z_\eta)
\otimes_\Z\Q$ a totally indefinite quaternion algebra over $\Q$.
\item[b.] for $d$ even, ${\rm End}(Z_\eta)
\otimes_\Z\Q$ a totally definite quaternion algebra over $\Q$.
\end{enumerate}
\end{enumerate}
If $g:Z\to Y$ is any family of Abelian varieties with a strictly maximal
Higgs field, then by \cite{VZ4} there exists some $A$ such that, replacing $Y$
by an \'etale covering, $Y=Y_A$ and $g:Z\to Y$ is isogenous to
$$
Z_A\times_{Y} \cdots \times_{Y}Z_A \>>> Y.
$$
In particular, the dimension $n$ of the fibres of $g$ is even.
\begin{example}\label{examplec}
For $n$ even, there exists a smooth family $f:X\to Y$ of $n$-dimensional manifolds
$F$ with $h^{n,0}(F)=1$ and $h^{k,0}\neq 0$ for all even $k$, such that the Higgs field
of $R^nf_*\Q_X$ is strictly maximal.
\end{example}
\begin{proof}
Starting from a smooth family
$$
g:Z=Z_A\times_{Y} \cdots \times_{Y}Z_A \>>> Y
$$
let $g':Z'\to Y$ be the family obtained by
blowing up the fixed points $\Gamma$ of the involution $\iota$ of $Z$ over $Y$.
We will assume for simplicity that those fixed points are the union of
the images of sections of $g$, a condition always satisfied after replacing $Y$
by an \'etale covering.

The involution $\iota$ acts on the relative Zariski tangent space of $\Gamma$
by multiplication with $-1$, hence it induces an action on $Z'$,
denoted by $\iota'$. The restriction of $\iota'$ to the exceptional
divisor $E$ is trivial, and $\iota'$ acts
fixed point free on $Z'\setminus E$. If $\V^k$ denotes the variation of Hodge
structures of weight $k$ of $Z\to Y$, one has
$$
R^kg'_*\C_X\otimes\sO_Y= R^k g'_*\Omega^\bullet_{Z'/Y}
= \V^k\otimes \sO_Y \oplus \T^k\otimes \sO_Y.$$
Here $\T^k$ is zero for $k$ odd, and it is a
trivial local system, concentrated in bidegree $(\frac{k}{2},\frac{k}{2})$ for
$k$ even. Moreover, $\iota'$ acts trivially on $\T^k$ and on $\V^k$ it acts
by multiplication with $(-1)^k$.

The quotient, $X=\Z'/\iota'$ is smooth over $Y$, and the ramification locus
of $\tau: Z'\to X$ is $E$. One finds
$$
R^k f_*\Omega^\bullet_{X/Y} = \V^k\otimes \sO_Y \oplus \T^k\otimes \sO_Y,
$$
for $k$ even and $=0$ for $k$ odd.
\end{proof}
The last example shows that for $n$ even
there exist smooth families of $n$-folds of Kodaira dimension zero, with
$h^{n,0}=1$ and with a weight $n$ variation of Hodge structures with a strictly
maximal Higgs field, whereas the same is excluded for $n$ odd by \ref{empty2}.
\begin{problem} Let $n>2$ be even.
Do there exist smooth families of Calabi-Yau $n$-folds over projective curves
with a weight $n$ variation of Hodge structures, with a strictly maximal Higgs field?
\end{problem}
\bibliographystyle{plain}

\end{document}